\documentclass[a4paper, leqno]{mathscan}

\newtheorem{thm}{Theorem}[section] 
       \newtheorem{pro}[thm]{Proposition} 
       \newtheorem{cor}[thm]{Corollary}    
       \newtheorem{lem}[thm]{Lemma}        

       \theoremstyle{definition} 
         
\def\proof{\noindent {\it Proof. }}

\renewcommand{\mathbb}{\mathsf}
\newcommand{\cal}{\mathcal}

\begin{document}

\title[MULTIVARIATE MEIXNER-POLLACZEK POLYNOMIALS]{HERMITIAN SYMMETRIC SPACES OF TUBE TYPE  AND MULTIVARIATE MEIXNER-POLLACZEK POLYNOMIALS}

\author{Jacques Faraut}  

\address{
\newline Institut de Math\'ematiques de Jussieu
\newline Universit\'e Pierre et Marie Curie
\newline 4 place Jussieu, case 247,
75 252 Paris cedex 05, 
\newline France}

\email{jacques.faraut@imj-prg.fr}

\author{Masato Wakayama}

\address{
\newline Institute of Mathematics for Industry
\newline Kyushu University
\newline Motooka, Nishi-ku, Fukuoka 819-0395, 
\newline Japan}

\email{wakayama@imi.kyushu-u.ac.jp}




\thanks{}

\begin{abstract} Harmonic analysis on Hermitian symmetric spaces of tube type is a natural framework for introducing multivariate Meixner-Pollaczek polynomials. Their main properties are established in this setting: orthogonality, generating and determinantal formulae, difference equations. For proving these properties we use the composition of the following transformations: Cayley transform, Laplace transform,
and spherical Fourier transform associated to Hermitian symmetric spaces of tube type.
In particular the difference equation for the multivariate Meixner-Pollaczek polynomials is obtained
from an Euler type equation on a bounded symmetric domain..
\end{abstract}

\maketitle

\centerline{\bf Contents}

1. Introduction

2. Spherical Fourier analysis on a symmetric cone

3. Multivariate Meixner-Pollaczek polynomials $Q_{\bf m}^{(\nu )}$

4. Multivariate Meixner-Pollaczek polynomials $Q_{\bf m}^{(\nu ,\theta)}$

5. Determinantal formulae

6. Difference equation for the multivariate Meixner-Pollaczek polynomials

7. The symmetries $S_{\nu }^{(i)}$ and the Hankel transform

8. Proof of the difference equation

9. Pieri's formula for the multivariate Meixner-Pollaczek polynomials


\section{Introduction}

The one variable Meixner-Pollaczek polynomials $P_m^{\alpha }(\lambda ;\phi)$ can be 
defined by the Gaussian hypergeometric representation as 
$$
P_m^{({\nu \over 2})}\bigl(\lambda ; \phi\bigr)
={(\nu )_m\over m!}e^{im\phi }
{}_2F_1\bigl(-m,{\nu \over 2}+i\lambda;\nu ;1-e^{-2i\phi }\bigr).
$$
For $\phi={\pi \over 2}$ the Meixner-Pollaczek polynomials 
$P_m^{({\nu \over 2})}\bigl(\lambda ; {\pi \over 2}\bigr)$  
are also obtained as Mellin transforms of Laguerre functions. Their main properties follow from this fact: hypergeometric representation above, orthogonality, generating formula, difference equation, and three terms relation (see \cite{AAR1999} p.348-349). 

These polynomials $P_m^{({\nu \over 2})}\bigl(\lambda ;{\pi \over 2}\bigr)$ have been generalized to the multivariate case. In fact, the multivariable Meixner-Pollaczek  (symmetric) polynomials
have been essentially considered in the setting of the Fourier analysis on Riemannian symmetric spaces in several papers: See Peetre-Zhang \cite{PZ1992} (Appendix 2: A class of hypergeometric orthogonal polynomials), \O rsted-Zhang \cite{OZ1994}, section 3.4, Zhang \cite{Z2002} and Davidson-\'Olafsson-Zhang \cite{DOZ2003}. Also, see the papers by Davidson-\'Olafsson \cite{DO2003} and Aristidou-David\-son-\'Olafsson \cite{ADO2006}. 
Further, for an arbitrary real value of the multiplicity $d$, the multivariate Meixner-Pollaczek polynomials are defined by  Sahi-Zhang \cite{SZ2007} in the setting of Heckman-Opdam and Cherednik-Opdam transforms, related to symmetric and non-symmetric Jack polynomials, and generating formulae for them are established. However the case where the parameter $\phi$ is involved has not been studied so far. Moreover, once we define the multivariate Meixner-Pollaczek polynomials with parameter $\phi$, it is also important to clarify a geometric meaning of the parameter. Establishing a natural setting for the study of multivariate Meixner-Pollaczek polynomials with such parameter, one can expect to obtain wider applications such as a study of multi-dimensional L\'evi-process, in particular, introducing multi-dimensional Meixner process (see \cite{S2000}
for the one dimensional case).  

The purpose of this article is to provide a geometric framework for introducing 
the multivariate Meixner-Pollaczek polynomials (with parameter $\phi$) and study their fundamental properties. Our analysis may explain much simpler geometric 
understanding of several basic properties of the multivariate Meixner-Pollaczek polynomials 
than ever, 
even in the case $\phi ={\pi \over 2}$. 
For instance, the $\mathfrak{S}_n$-invariant difference operator of which the multivariate Meixner-Pollaczek polynomials are eigenfunctions can be 
understood by an image of the Euler operator under the composition of three intertwiners: the Cayley transform, the Laplace transform and the spherical Fourier transform. 
In particular, the multivariate Meixner-Pollaczek polynomials are spherical Fourier transforms
of multivariate Laguerre functions.

In Section~2 we recall the basic facts about the spherical Fourier analysis on a symmetric cone. In Section~3 we define the multivariate Meixner-Pollaczek polynomials $Q_{\bf m}^{(\nu )}({\bf s})$ (the case $\phi ={\pi \over 2}$), where $\bf m$ is a partition, prove that they are orthogonal with respect to a measure $M_{\nu }$ on $\mathbb{R}^n$, and establish a generating formula.

In Section~4, adding a real parameter $\theta $ (instead of $\phi =\theta +{\pi \over 2}$), we introduce the symmetric polynomials $Q_{\bf m}^{(\nu ,\theta )}({\bf s})$ in the variables ${\bf s}=(s_1,\ldots,s_n)$ 
($Q_{\bf m}^{(\nu )}=Q_{\bf m}^{(\nu ,0)}$). In the one variable case
$$q_m^{(\nu ,\theta )}(s )=(-i)^mP_m^{({\nu \over 2})}\bigl(-is ; \theta +{\pi \over 2}\bigr).$$
The orthogonality property for the polynomials $Q_{\bf m}^{(\nu ,\theta )}({\bf s})$ is obtained by using
a Gutzmer formula for the spherical Fourier transform. A generating formula is obtained for these polynomials. 
In case of the multiplicity $d=2$, we establish in Section~5 determinantal formulae for multivariate Laguerre and Meixner-Pollaczek polynomials. Sections~6, 7, and 8 are devoted to a difference equation satisfied by the polynomials $Q_{\bf m}^{(\nu ,\theta )}({\bf s})$.
Starting from an Euler-type equation involving the parameter $\theta $, this difference equation is obtained in three steps, corresponding to a Cayley transform, an inverse Laplace transform, and a spherical Fourier transform for symmetric cones. The symmetry $\theta \mapsto -\theta$ in the parameter is related to geometric symmetries and to 
a generalized Tricomi theorem for the Hankel transform on a symmetric cone. In the last section we show that multivariate Meixner-Pollaczek polynomials
satisfy a Pieri's formula. In the one variable case it reduces to the three terms relation satisfied
by the classical Meixner-Pollacek polynomials. 
\section{Spherical Fourier analysis on a symmetric cone}
A reference for this preliminary section is \cite{FK1994}.
We consider an irreducible symmetric cone $\Omega $ in a Euclidean Jordan algebra $V$.
We denote by $G$ the identity component in the group $G(\Omega )$ of linear automorphisms of $\Omega$, and $K\subset G$ is the isotropy subgroup of the unit element $e\in V$.

The Gindikin gamma function $\Gamma _{\Omega }$ of the cone $\Omega $ 
will be the cornerstone of the analysis we will developp.
It is defined, for ${\bf s}\in \mathbb{C}^n$, with ${\rm Re}\, s_j>{d\over 2}(j-1)$, by
$$\Gamma _{\Omega } ({\bf s})=\int _{\Omega }e^{-{\rm tr}\, (u)}\Delta _{\bf s}(u)\Delta (u)^{-{N\over n}}m(du).$$ 
The notation ${\rm tr}\, (u)$ and $\Delta (u)$ denote the trace and the determinant with respect to the Jordan algebra structure, $\Delta _{\bf s}$ is the power function, $N$ and $n$ are the dimension and the rank
of $V$, and $m$ is the Euclidean measure associated to the Euclidean structure on $V$ given by
$(u|v)={\rm tr}\, (uv)$. Its evaluation gives
$$\Gamma _{\Omega }({\bf s})
=(2\pi )^{N-n\over 2}\prod _{j=1}^n \Gamma \bigl(s_j-{d\over 2}(j-1)\bigr),$$
where $d$ is the multiplicity, related to $N$ and $n$ by the relation $N=n+{d\over 2}n(n-1)$.
The spherical function $\varphi _{\bf s}$, for ${\bf s}\in \mathbb{C}^n$, is defined on $\Omega $
by
$$\varphi _{\bf s}(u)=\int _K\Delta _{{\bf s}+\rho }(k\cdot u)dk,$$
where $\rho =(\rho _1,\ldots ,\rho _n)$, $\rho _j={d\over 4}(2j-n-1)$, and $dk$ is the normalized Haar measure on the compact group $K$.
The algebra $\mathbb{D}(\Omega )$ of $G$-invariant differential operators on $\Omega $ is commutative, and the spherical function $\varphi _{\bf s}$ is an eigenfunction of every 
$D\in \mathbb{D}(\Omega )$: 
$$D\varphi _{\bf s}=\gamma _{D}({\bf s})\varphi _{\bf s}.$$
The function $\gamma _D$ is a symmetric polynomial function, and the map
$D\mapsto \gamma _D$ is an algebra isomorphism from $\mathbb{D}(\Omega )$ onto 
the algebra ${\cal P}(\mathbb{C}^n)^{\mathfrak{S}_n}$ of symmetric polynomial functions,
a special case of the Harish-Chandra isomorphism.
The spherical Fourier transform ${\cal F}\psi $ of a $K$-invariant function $\psi $ on $\Omega $ is given by
$${\cal F}\psi ({\bf s})=\int _{\Omega } \psi (u)\varphi _{\bf s}(u)
\Delta ^{-{N\over n}}(u)m(du).$$
Hence, for $\psi (u)=e^{-{\rm tr}\,  u}\Delta ^{\nu \over 2}(u)$ ($\nu >{d\over 2}(n-1)$), then
$${\cal F}\psi ({\bf s})=\Gamma _{\Omega } ({\bf s}+{\nu \over 2}+\rho \bigr)
=(2\pi )^{N-n\over 2}\prod _{j=1}^n \Gamma \bigl(s_j+{\nu \over 2}-{d\over 4}(n-1)\bigr).$$
For an invariant differential operator $D\in \mathbb{D}(\Omega )$,
${\cal F}(D\psi )=\gamma _D(-{\bf s}){\cal F}\psi $.
The space ${\cal P}(V)$ of polynomials on $V$ decomposes multiplicity free under $G$ as
$${\cal P}(V)=\bigoplus _{\bf m}{\cal P}_{\bf m},$$
where ${\cal P}_{\bf m}$ is a finite dimensional subspace, irreducible under $G$. The parameter $\bf m$
is a partition: ${\bf m}=(m_1,\ldots ,m_n)\in \mathbb{N}^n$, $m_1\geq \cdots \geq m_n$. The polynomials in ${\cal P}_{\bf m}$ are homogeneous of degree
$|{\bf m}|:=m_1+\cdots +m_n$. The subspace
${\cal P}_{\bf m}^K$ of $K$-invariant polynomials in ${\cal P}_{\bf m}$ is one dimensional, generated by the spherical polynomial $\Phi _{\bf m}$, normalized by the condition $\Phi _{\bf m}(e)=1$,
and $\Phi _{\bf m}=\varphi _{{\bf m}-\rho }$.
There is a unique invariant differential operator $D^{\bf m}$ such that
$$D^{\bf m}\psi (e)=\Bigl(\Phi _{\bf m}\bigl({\partial \over \partial u}\bigr)\psi  \Bigr)(e).$$
We will write $\gamma _{\bf m}=\gamma _{D^{\bf m}}$.
Observe that, for $n=1$, $\Phi _m(u)=u^m$, and
$$D^m=u^m\Bigl({d\over du}\Bigr)^m,\ \gamma _m(s)=[s]_m:=s(s-1)\ldots (s-m+1).$$
The classical Pochhammer symbol $(\alpha )_m:=\alpha (\alpha +1)\ldots (\alpha +m-1)$
generalizes as follows: for $\alpha \in \mathbb{C}$ and a partition $\mathbf{m}$,
$$(\alpha )_{\bf m}={\Gamma _{\Omega }({\bf m}+\alpha )\over \Gamma _{\Omega } (\alpha )}
=\prod _{i=1}^n\Bigl(\alpha -(i-1){d\over 2}\Bigr)_{m_i}.$$
If a $K$-invariant function $\psi $ is analytic in a neighborhood of $e$,
it admits a spherical Taylor expansion near $e$:
$$\psi (e+v)=\sum _{\bf m} d_{\bf m}{1\over \bigl({N\over n}\bigr)_{\bf m}}
D^{\bf m}\psi (e)\Phi _{\bf m} (v),$$
where $d_{\bf m}$ is the dimension of ${\cal P}_{\bf m}$.
In particular, for $\psi =\varphi _{\bf s}$, a spherical function, 
$$\varphi _{\bf s}(e+v)=\sum _{\bf m} d_{\bf m} {1\over \bigl({N\over n}\bigr)_{\bf m}}
\gamma _{\bf m}({\bf s})\Phi _{\bf m}(v).$$
For $\psi =\Phi _{\bf m}=\varphi _{{\bf m}-\rho }$, we get the spherical binomial formula
$$\Phi _{\bf m}(e+v)=\sum _{{\bf k}\subset {\bf m}}{{\bf m}\choose {\bf k}}\Phi _{\bf k} (v).$$
In fact the generalized binomial coefficient 
$${{\bf m}\choose {\bf k}}={d_{\bf k}}{1\over \bigl({N\over n}\bigr)_{\bf k}}
\gamma _{\bf k}({\bf m}-\rho )$$
vanishes if ${\bf k}\not\subset {\bf m}$.
\section{Multivariate Meixner-Pollaczek polynomials $Q_{\bf m}^{(\nu )}$}
For $n=1$, we define the Meixner-Pollaczek polynomial $q_m^{(\nu )}$ as follows:
$$q_m^{(\nu )}(s)={(\nu )_m\over m!}{} _2F_1(-m,s+{\nu \over 2};\nu ;2).$$
This definition slightly differs from the classical one $P_m^{\alpha }(\lambda ;\phi )$:
$$q_m^{(\nu )}(i\lambda )=(-i)^mP_m^{\nu \over 2}\bigl(\lambda ;{\pi \over 2}\bigr).$$
(see for instance \cite{AAR1999}, p.348.)
Its expansion can be written
$$q_m^{(\nu )}(s)={(\nu )_m\over m!}
\sum _{k=0}^m {[m]_k\big[-s-{\nu \over 2}\big]_k\over (\nu )_k}{1\over k!}2^k.$$
The polynomials $q_m^{(\nu )}(i\lambda )$ are orthogonal with respect to the weight on $\mathbb{R}$
$$|\Gamma \bigl(i\lambda +{\nu \over 2}\bigr)|^2\quad (\nu >0).$$
We define the multivariate Meixner-Pollaczek polynomial $Q_{\bf m}^{(\nu )}$ as
 the following symmetric polynomial in $n$ variables: 
 $$Q_{\bf m}^{(\nu )}({\bf s})={(\nu )_{\bf m}\over \bigl({N\over n}\bigr)_{\bf m}}
 \sum _{{\bf k}\subset {\bf m}} d_{\bf k}
 {\gamma _{\bf k}({\bf m}-\rho )\gamma _{\bf k}\bigl(-{\bf s}-{\nu \over 2}\bigr)\over (\nu )_{\bf k}}
 {1\over \bigl({N\over n}\bigr)_{\bf k}}2^{|{\bf k}|}.$$
 For $\nu >{d\over 2}(n-1)$ let us denote by $M_{\nu }(d\lambda )$ the probability measure on $\mathbb{R}^n$ given by
 $$M_{\nu }(d\lambda )=
 {1\over Z_{\nu }} \prod _{j=1}^n \Bigl|\Gamma \bigl(i\lambda _j+{\nu \over 2}-{d\over 4}(n-1)\bigr)\Bigr|^2
 {1\over |c(i\lambda )|^2}m(d\lambda ),$$
 where 
 $$Z_{\nu }=\int _{\mathbb{R}^n}
 \prod _{j=1}^n \Bigl|\Gamma \bigl(i\lambda _j+{\nu \over 2}-{d\over 4}(n-1)\bigr)\Bigr|^2
 {1\over |c(i\lambda )|^2}m(d\lambda ),$$
 and $c$ is the Harish-Chandra function for the symmetric cone $\Omega $,
 $$c({\bf s})=c_0\prod _{j<k}B\Bigl(s_j-s_k,{d\over 2}\Bigr).$$
 ($B$ is the Euler beta function, the constant $c_0$ is such that $c(-\rho )=1$, see Section XIV.5
 in \cite{FK1994}.)
 The constant $Z_{\nu }$ can be evaluated by using the spherical Plancherel formula, applied to the function 
$\psi (u)=e^{-{\rm tr}\, u}\Delta (u)^{\nu \over 2}$: 
\begin{eqnarray*}
& &\int _{\Omega } e^{-2{\rm tr}\,  u} \Delta (u)^{\nu -{N\over n}}m(du)\\
& &=(2\pi )^{N-2n}\int _{\mathbb{R}^n} 
 \prod _{j=1}^n |\Gamma (i\lambda _j+{\nu \over 2}-{d\over 4}(n-1)|^2
 {1\over |c(i\lambda )|^2}m(d\lambda ).
\end{eqnarray*}
 Therefore
 $$Z_{\nu }=(2\pi )^{2n-N} 2^{-n\nu }\Gamma _{\Omega }(\nu ).$$
Next statement involves the geometry of the Hermitian symmetric space of tube type associated to the symmetric cone $\Omega $. The map $z\mapsto (z-e)(z+e)^{-1}$ maps the tube domain
$T_{\Omega }=\Omega +iV\subset V_{\mathbb{C}}$ onto the bounded Hermitian symmetric domain ${\cal D}$. Its inverse is the Cayley transform:
$$c(w)=(e+w)(e-w)^{-1}.$$
 
\begin{thm}\label{MPOG1}
 Assume $\nu >{d\over 2}(n-1)$.
 
 {\rm (i)} The multivariate Meixner-Pollaczek polynomials $Q_{\bf m}^{(\nu )}(i\lambda )$ form an orthogonal basis
 of $L^2(\mathbb{R}^n,M_{\nu })^{\mathfrak{S}_n}$.
 The norm of $Q_{\bf m}^{(\nu )}$ is given by:
 $$\int _{\mathbb{R}^n}|Q_{\bf m}^{(\nu )}(i\lambda )|^2M_{\nu }(d\lambda )
 ={1\over d_{\bf m}}{(\nu )_{\bf m}\over \bigl({N\over n}\bigr)_{\bf m}}.$$
 
 {\rm (ii)} The polynomials $Q_{\bf m}^{(\nu )}$ admit the following generating formula:
 for ${\bf s}\in \mathbb{C}^n$, $w\in {\cal D}$,
 $$\sum _{\bf m}d_{\bf m}Q_{\bf m}^{(\nu )}({\bf s})\Phi _{\bf m} (w)
 =\Delta (e-w^2)^{-{\nu \over 2}}\varphi _{\bf s}\bigl(c(w)^{-1}\bigr).$$
 \end{thm} 
 
\proof
a) For $\nu >2{N\over n}-1=1+d(n-1)$, ${\cal H}_{\nu }^2({\cal D})$ denotes the weighted Bergman space of holomorphic functions $f$ on $\cal D$ such that
$$\|f\|_{\nu }^2:=a_{\nu }^{(1)} \int _{\cal D}|f(w)|^2h(w)^{\nu -2{N\over n}}m(dw)<\infty .$$
The constant 
$$a_{\nu }^{(1)}={1\over \pi ^n}
{\Gamma _{\Omega }(\nu )\over \Gamma _{\Omega }\bigl(\nu -{N\over n}\bigr)}$$ 
is such that the function $\Phi _0\equiv 1$ has norm 1.
Recall that $h(w)=h(w,w)$, where $h('w,w')$ is a polynomial holomorphic in $w$, antiholomorphic in $w'$, such that, for $w$ invertible, $h(w,w')=\Delta (w)\Delta (w^{-1}-\bar w')$
($\bar w'$ is the complex conjugate of $w'$ with respect to the real form $V$ of $V_{\mathbb{C}}$).
The spherical polynomials $\Phi _{\bf m}$ form an orthogonal basis of the space 
${\cal H}_{\nu }^2({\cal D})^K$ of $K$-invariant functions in ${\cal H}_{\nu }^2({\cal D})$,
and
$$\| \Phi _{\bf m}\|_{\nu }^2={1\over d_{\bf m}}{\bigl({N\over n}\bigr)_{\bf m}\over (\nu )_{\bf m}}.
\eqno (3.1)$$
The reproducing kernel of ${\cal H}_{\nu }^2({\cal D})$ is given by
${\cal K}_{\nu } (w,w')=h(w,w')^{-\nu }$.
By an integration over $K$ one obtains:
$${\cal G}_{\nu }^{(1)} (\zeta ,w)
:=\sum _{\bf m}d_{\bf m}{(\nu )_{\bf m}\over \bigl({N\over n}\bigr)_{\bf m }}
\Phi _{\bf m}(\zeta )\Phi _{\bf m}(w)
=\int _K h(w,k\bar \zeta )^{-\nu }dk 
.\eqno (3.2)$$

b) For a function $f$ holomorphic in $\cal D$, one defines the function $F=C_{\nu } f$ on $T_{\Omega }$
by
$$F(z)=\bigl(C_{\nu } f)(z)=\Delta \bigl({z+e\over 2}\bigr)^{-\nu } f\bigl((z-e)(z+e)^{-1}\bigr).$$
The map $C_{\nu }$ is a unitary isomorphism from ${\cal H}_{\nu }^2({\cal D})$ onto the space ${\cal H}_{\nu }^2(T_{\Omega })$ of holomorphic functions on $T_{\Omega }$ such that
$$\|F\|_{\nu }^2 
:=a_{\nu }^{(2)} \int _{T_{\Omega }}|F(z)|^2\Delta (x)^{\nu -2{N\over n}}m(dz)<\infty .$$
The constant
$$a_{\nu }^{(2)}={1\over (4\pi )^n}
{\Gamma _{\Omega }(\nu )\over \Gamma _{\Omega } \bigl(\nu -{N\over n}\bigr)},$$
is such that the function 
$$F_0^{(\nu )}=C_{\nu } \Phi _0,\quad F_0^{(\nu )}(z)=\Delta \bigl({z+e\over 2}\bigr)^{-\nu },$$
has norm 1. The functions $F_{\bf m}^{(\nu )}=C_{\nu }\Phi _{\bf m}$ form an orthogonal basis of 
the space ${\cal H}_{\nu }^2(T_{\Omega })^K$ of $K$-invariant functions in ${\cal H}_{\nu }^2(T_{\Omega })$,
and it follows from (3.1) that
$$\|F_{\bf m}^{ (\nu )}\|_{\nu }^2={1\over d_{\bf m}}{\bigl({N\over n}\bigr)_{\bf m}\over (\nu )_{\bf m}}.
\eqno (3.3)$$
Performing the transform $C_{\nu }$ with respect to $\zeta $ in (3.2) we get
a generating formula for the functions $F_{\bf m}^{(\nu )}$: for $w\in {\cal D}$,
$z\in T_{\Omega }$,
\begin{eqnarray*}
{\cal G}_{\nu }^{(2)} (z,w)
&:=&\sum _{\bf m}d_{\bf m} {(\nu )_{\bf m}\over \bigl({N\over n}\bigr)_{\bf m}}
\Phi _{\bf m}(w)F_{\bf m}^{(\nu )}(z) \\
&=&\Delta \bigl({e-w\over 2}\bigr)^{-\nu }
\int _K \Delta \bigl(k\cdot z+c(w)\bigr)^{-\nu }dk. \qquad \qquad (3.4) 
\end{eqnarray*}

c) The functions in ${\cal H}_{\nu }^2(T_{\Omega })$ admit a Laplace integral representation. The modified Laplace transform ${\cal L}_{\nu }$, given, for a function $\psi $ on $\Omega $, by
$$({\cal L}_{\nu }) \psi (z)
=a_{\nu }^{(3)} \int _{\Omega } e^{(z|u)}\psi (u)\Delta (u)^{\nu -{N\over n}}m(du),$$
is an isometric isomorphism from the space $L_{\nu }^2(\Omega )$ of measurable functions $\psi $ on $\Omega $ such that
$$\|\psi \|_{\nu }^2:=a_{\nu }^{(3)} \int _{\Omega } |\psi (u)|^2 \Delta (u)^{\nu -{N\over n}}m(du)<\infty ,$$
onto ${\cal H}_{\nu }^2(T_{\Omega })$. The constant
$a_{\nu }^{(3)}=2^{n\nu }/\Gamma _{\Omega }(\nu )$
is such that the function $\Psi _0(u)=e^{-{\rm tr}\, u}$ has norm 1, and then ${\cal L}_{\nu } \Psi _0=F_0$.
By the binomial formula
\begin{eqnarray*}
F_{\bf m}^{(\nu )}(z)
&=&\Delta \bigl({z+e\over 2}\bigr)^{-\nu }
\Phi _{\bf m}\bigl((z-e)(z+e)^{-1}\bigr)
=\Delta \bigl({z+e\over 2}\bigr)^{-\nu }
\Phi _{\bf m}\bigl(e-2(z+e)^{-1}\bigr)\\
&=&\sum _{{\bf k}\subset {\bf m}}(-1)^{|{\bf k}|}{{\bf m}\choose {\bf k}}
\Phi _{\bf k}\bigl(2(z+e)^{-1}\bigr)\Delta \bigl(2(e+z)^{-1}\bigr)^{\nu }.
\end{eqnarray*}
By Lemma XI.2.3 in \cite{FK1994} we have the following 

\begin{lem}\label{LaplaceImage}
$${\cal L}_{\nu } \bigl(e^{-\rm{tr}\, u}\Phi _{\bf m}\bigr)(z)
=(\nu )_{\bf m} \Phi _{\bf m}\bigl((z+e)^{-1}\bigr)\Delta \bigl(2(e+z)^{-1}\bigr)^{\nu }.$$
\end{lem}

By Lemma \ref{LaplaceImage} the function
$$\Psi _{\bf m}^{(\nu )}
={(\nu )_{\bf m}\over \bigl({N\over n}\bigr)_{\bf m}}{\cal L}_{\nu }^{-1}\bigl(F_{\bf m}^{(\nu )}\bigr).$$
is the Laguerre function given by
$$\Psi _{\bf m}^{(\nu )}(u)=e^{-{\rm tr}\, u}L_{\bf m}^{(\nu -1)}(2u),$$
where $L_{\bf m}^{(\nu -1)}$ is the multivariate Laguerre polynomial 

\begin{eqnarray*}
L_{\bf m}^{(\nu -1)}(x)
&=&{(\nu )_{\bf m}\over \bigl({N\over n}\bigr)_{\bf m}}
\sum _{{\bf k}\subset {\bf m}}{{\bf m}\choose {\bf k}}{1\over (\nu )_{\bf k}}\Phi _{\bf k}(-x)\\
&=&{(\nu )_{\bf m}\over \bigl({N\over n}\bigr)_{\bf m}}
\sum _{{\bf k}\subset {\bf m}} 
d_{\bf k}{\gamma _{\bf k}({\bf m}-\rho )\over (\nu )_{\bf k}}{1\over \bigl({N\over n}\bigr)_{\bf k}}
\Phi _{\bf k}(-x).
\end{eqnarray*}

\begin{pro}\label{LOG}
{\rm (i)} The multivariate Laguerre functions $\Psi _{\bf m}^{(\nu )}$ form an orthogonal basis of $L_{\nu }^2(\Omega )^K$, and
$$\|\Psi _{\bf m}^{(\nu )}\|_{\nu }^2={1\over d_{\bf m}}{(\nu )_{\bf m}\over \bigl({N\over n}\bigr)_{\bf m}}.
\eqno (3.5)$$
{\rm (ii)} The fonctions $\Psi _{\bf m}^{(\nu )}$ admit the following generating formula:
for $u\in \Omega $, $w\in {\cal D}$,
$${\cal G}_{\nu }^{(3)}(u,w)
:=\sum _{\bf m} d_{\bf m}\Psi _{\bf m}^{(\nu )}(u)\Phi _{\bf m}(w)
=\Delta (e-w)^{-\nu }\int _K e^{-\bigl(k\cdot u|c(w)\bigr)}dk
.\eqno (3.6)$$
The generating formula can also be written
$$\Delta (e-w)^{-\nu } \int _K e^{(k\cdot x|w(e-w)^{-1})}dk
=\sum _{\bf m} d_{\bf m}L_{\bf m}^{(\nu -1)}(x)\Phi _{\bf m}(w).\eqno (3.6')$$
\end{pro}

Formula (3,6') is proposed  as an exercise in \cite{FK1994} (Exercise 3, p.347). It is a special case of formula (4.4) in \cite{BF1997}.
\smallskip

\proof
Part (i) follows from the fact that ${\cal L}_{\nu }$ is a unitary isomorphism from $L_{\nu }^2(\Omega )$
onto ${\cal H}_{\nu }^2(T_{\Omega })$, and from (3.3).

The modified Laplace transform of ${\cal G}_{\nu }^{(3)}(u,w)$ with respect to $u$ 
is equal to ${\cal G}_{\nu }^{(2)} (z,w)$,
and one gets (ii) from (2.4).
\qed

\medskip

d) We will evaluate the spherical Fourier transform of the Laguerre functions $\Psi _{\bf m}^{(\nu )}$.
We introduce now the modified spherical Fourier transform ${\cal F}_{\nu }$ as follows: for a function $\psi $ on $\Omega $,
$$({\cal F}_{\nu }\psi )({\bf s})
={1\over \Gamma _{\Omega }\bigl({\bf s}+{\nu \over 2}+\rho \bigr)}
\int _{\Omega } \psi (u)\varphi _{\bf s}(u)\Delta (u)^{{\nu \over 2}-{N\over n}}m(du).$$
Observe that ${\cal F}_{\nu } \Psi _0\equiv 1$.

\begin{lem}\label{Fourier}
For ${\rm Re}\, s_j>{d\over 4}(n-1)-{\nu \over 2}$,
$${\cal F}_{\nu }\bigl(e^{-{\rm tr}\,  u}\Phi _{\bf m}\bigr)({\bf s})
=(-1)^{|{\bf m}|}\gamma _{\bf m}\bigl(-{\bf s}-{\nu \over 2}\bigr)
.$$
\end{lem}

\proof
Let $\sigma _D(u,\xi )$ be the symbol of $D\in \mathbb{D}(\Omega )$, and $p(\xi )=\sigma _D(e,\xi )$ (See \cite{FK1994}, p.290). By the invariance property of $\sigma _D$, we have $\sigma _D(u,-e)=p(-u)$, and therefore
$De^{-{\rm tr}\, u}=p(-\xi )e^{-{\rm tr}\, u}$.
Hence, for $p(\xi )=\Phi _{\bf m}(\xi )$,
\begin{eqnarray*}
&{\cal F}_{\nu } (e^{-{\rm tr}\, u}\Phi _{\bf m})(s)
=(-1)^{|{\bf m}|}{\cal F}_{\nu } (D^{\bf m}e^{-{\rm tr}\, u})(s) \\
&=(-1)^{|{\bf m}|}\gamma _{\bf m}\bigl(-{\bf s}-{\nu \over 2}\bigr){\cal F}_{\nu } (e^{-{\rm tr}\, u})
=(-1)^{|{\bf m}|}\gamma _{\bf m}\bigl(-{\bf s}-{\nu \over 2}\bigr).
\end{eqnarray*}
\qed

From  Lemma \ref{Fourier} we obtain the evaluation of the spherical Fourier transform of the Laguerre functions:
For ${\rm Re}\, s_j>{d\over 4}(n-1)-{\nu \over 2}$,
$${\cal F}_{\nu } (\Psi _{\bf m}^{\nu })({\bf s})=Q_{\bf m}^{(\nu )}({\bf s})
.$$
By the spherical Plancherel formula and part (i) in Proposition \ref{LOG}, this proves part (i) of Theorem \ref{MPOG1},
for $\nu >1+d(n-1)$:
$$\int _{\mathbb{R}^n}|Q_{\bf m}^{(\nu )}(i\lambda )|^2M_{\nu }(d\lambda )
 ={1\over d_{\bf m}}{(\nu )_{\bf m}\over \bigl({N\over n}\bigr)_{\bf m}}.\eqno (3.7)$$
 By analytic continuation it holds for $\nu >{d\over 2}(n-1)$.
For proving part (ii) of Theorem 2.1 one performs the spherical Fourier transform to both handsides
of part (ii) in Proposition \ref{LOG}:
$${\cal G}_{\nu }^{(4)}({\bf s},w):=\sum _{\bf m}d_{\bf m}Q_{\bf m}^{(\nu )}({\bf s})\Phi _{\bf m} (w)
 =\Delta (e-w^2)^{-{\nu \over 2}}\varphi _{\bf s}\bigl(c(w)^{-1}\bigr).\eqno (3.8)$$
This finishes the proof of Theorem \ref{MPOG1}.
We remark that, 
in \cite{DOZ2003}, a different notation is used for the Meixner-Pollaczek polynomials: their polynomials $p_{\nu ,{\bf m}}$  (p. 179) are defined through the 
generating formula above and  
$p_{\nu ,{\bf m}}(i{\bf s})=d_{\bf m}Q_{\bf m}^{(\nu )}({\bf s})$.
\section{Multivariate Meixner-Pollaczek polynomials $Q_{\bf m}^{(\nu,\theta )}$}
The Meixner-Pollaczek polynomials $q_m^{(\nu )}$ we have considered at the beginning of Section~3 correspond to the special value $\phi={\pi \over 2}$ with the classical notation. Using instead
$\theta=\phi -{\pi \over 2}$, the more general one variable Meixner-Pollaczek polynomials can be written
\begin{eqnarray*}
q_m^{(\nu ,\theta )}(s)
&=&e^{im\theta }{(\nu )_m\over m!}
{}_2F_1(-m,s+{\nu \over 2};\nu;2e^{-i\theta }\cos \theta )\\
&=&e^{im\theta }{(\nu )_m\over m!}
\sum _{k=0}^m {[m]_k\big[-s-{\nu \over 2}\big]_k\over (\nu )_k}
{1\over k!}(2e^{-i\theta }\cos \theta )^k.
\end{eqnarray*}
In terms of the classical notation $P_m^{\alpha }(\lambda ;\phi )$
$$q_m^{(\nu,\theta )}(i\lambda )=(-i)^mP_m^{\nu \over 2}\bigl(\lambda ;\theta +{\pi \over 2}\bigr).$$
For $\nu >0$,  $|\theta |<{\pi \over 2}$,
the polynomials $q_m^{(\nu,\theta )}(i\lambda )$ are orthogonal with respect to the weight
$$e^{2\theta \lambda } \big|\Gamma \bigl(i\lambda +{\nu \over 2}\bigr)\big|^2.$$
In this section we consider the multivariate Meixner-Pollaczek polynomials $Q_{\bf m}^{(\nu,\theta )}$
defined by
$$Q_{\bf m}^{\nu,\theta )}({\bf s})= e^{i|{\bf m}|\theta }{(\nu )_{\bf m}\over \bigl({N\over n}\bigr)_{\bf m}}
 \sum _{{\bf k}\subset {\bf m}} d_{\bf k}
 {\gamma _{\bf k}({\bf m}-\rho )\gamma _{\bf k}\bigl(-{\bf s}-{\nu \over 2}\bigr)\over (\nu )_{\bf k}}
 {1\over \bigl({N\over n}\bigr)_{\bf k}}(2e^{-i\theta }\cos \theta )^{|{\bf k}|}.$$

 \begin{thm}\label{MPOG2}
 Assume $\nu >{d\over 2}(n-1)$, $|\theta |<{\pi \over 2}$.
 
 {\rm (i)} The multivariate Meixner-Pollaczek polynomials $Q_{\bf m}^{(\nu ,\theta )}(i\lambda )$ form an orthogonal basis
 of $L^2(\mathbb{R}^n,e^{2\theta (\lambda _1+\cdots +\lambda _n)}M_{\nu })^{\mathfrak{S}_n}$.
 The norm of $Q_{\bf m}^{(\nu ,\theta )}$ is given by:
 $$\int _{\mathbb{R}^n}|
 Q_{\bf m}^{(\nu ,\theta )}(i\lambda )|^2e^{2\theta (\lambda _1+\cdots +\lambda _n)}M_{\nu }(d\lambda )
 =(\cos \theta )^{-n\nu }{1\over d_{\bf m}}{(\nu )_{\bf m}\over \bigl({N\over n}\bigr)_{\bf m}}.$$
 
 {\rm (ii)} The polynomials $Q_{\bf m}^{(\nu ,\theta )}$ admit the following generating formula:
 for ${\bf s}\in \mathbb{C}^n$, $w\in {\cal D}$,
 $$\sum _{\bf m}d_{\bf m}Q_{\bf m}^{(\nu ,\theta )}({\bf s})\Phi _{\bf m} (w)
 =\Delta \bigl((e-e^{i\theta }w)(e+e^{-i\theta }w)\bigr)^{-{\nu \over 2}}
 \varphi _{\bf s}\bigl(c_{\theta }(w)^{-1}\bigr),$$
 where $c_{\theta }$ is the modified Cayley transform:
 $$c_{\theta }(w)=(e+e^{-i\theta }w)(e-e^{i\theta }w)^{-1}.$$
\end{thm}
 
We will prove Theorem \ref{MPOG2} in several steps.
 
 a) Let us define the Laguerre functions $\Psi _{\bf m}^{(\nu ,\theta )}$:
 $$\Psi _{\bf m}^{(\nu ,\theta )}(u)=e^{i|{\bf m}|\theta }
 e^{-{\rm tr}\, u}L_{\bf m}^{(\nu -1)}(2e^{-i\theta }\cos \theta  \, u).$$
 
For  functions $\psi $ on $V$ of the form
$\psi (u)=e^{-{\rm tr}\, u}p(u)$, where $p$ is a polynomial, define the inner product
$$(\psi _1|\psi _2)_{(\nu ,\theta )}={2^{n\nu }\over \Gamma _{\Omega }(\nu )}
\int _{\Omega }\psi _1(e^{i\theta }u)\overline{\psi _2(e^{i\theta } u)}\Delta (u)^{\nu -{N\over n}}m(du).$$

\begin{pro}\label{LOG2}
{\rm (i)} The Laguerre functions $\Psi _{\bf m}^{(\nu, \theta )}$ are orthogonal with respect to the 
inner product $(\cdot |\cdot )_{(\nu ,\theta )}$. Furthermore
$$\|\Psi _{\bf m}^{(\nu ,\theta )}\|_{(\nu ,\theta )}^2
=(\cos \theta )^{-n\nu }{1\over d_{\bf m}}{(\nu )_{\bf m}\over \bigl({N\over n}\bigr)_{\bf m}}.$$

{\rm (ii)} The Laguerre functions $\Psi _{\bf m}^{(\nu ,\theta )}$ satisfy the following generating formula:
for $u\in \Omega$, $w\in {\cal D}$,
$${\cal G}_{\nu ,\theta }^{(3)}(u,w):=
\sum _{\bf m} d_{\bf m} \Psi _{\bf m}^{(\nu,\theta )}(u)\Phi _{\bf m}(w)
=\Delta (e-e^{i\theta }w)^{-\nu }\int _Ke^{\bigl(k\cdot u|c_{\theta }(w)\bigr)}dk
.$$
\end{pro}

\proof
(i) Put $\alpha =e^{i\theta }$, $\beta =2e^{-i\theta }\cos \theta$. For two polynomials $p_1$ and $p_2$ consider the functions
$$\psi _1^{(\theta )}(u)=e^{-{\rm tr}\, u}p_1(\beta u),\ \psi _2^{(\theta )}(u)=e^{-{\rm tr}\, u}p _2(\beta u),$$
and their inner product
$$(\psi _1^{(\theta )}|\psi _2^{(\theta )})_{\nu ,\theta }
={2^{n\nu }\over \Gamma _{\Omega } (\nu )}
\int _{\Omega } e^{-\alpha {\rm tr}\, u}p_1(\beta \alpha u)\overline{e^{-\alpha {\rm tr}\, u} p_2(\beta \alpha u)}
\Delta (u)^{\nu -{N\over n}}m(du).$$
Observe that $\beta \alpha =2\cos \theta $, $\alpha +\bar \alpha =2\cos \theta $. Hence

\begin{eqnarray*}
(\psi _1^{(\theta )}|\psi _2^{(\theta )})_{\nu ,\theta }
&=&{2^{n\nu }\over \Gamma_{\Omega }(\nu )}
\int _{\Omega } e^{-2\cos \theta {\rm tr}\, u}p_1(2\cos \theta u)\overline{p_2(2\cos \theta u)}
\Delta (u)^{\nu -{n\over N}}m(du)\\
&=&{2^{n\nu }\over \Gamma _{\Omega } (\nu )}(\cos \theta )^{-n\nu }
\int _{\Omega }e^{-2{\rm tr}\, v}p_1(2v)\overline{p_2(2v)}\Delta (v)^{\nu -{N\over n}}m(dv)\\
&=&(\cos \theta )^{-n\nu }(\psi _1^{(0)}|\psi _2^{(0)}).
\end{eqnarray*}

Take
$$p_1(u)=L_{\bf p}^{(\nu -1)}(u),\ p_2(u)=L_{\bf q}^{(\nu -1)}(u).$$
Then, by part (i) of Proposition \ref{LOG}, the statement (i) is proven.
 
(ii) The sum in the generating formula can be written
$$\sum _{\bf m} d_{\bf m} e^{-{\rm tr}\, u}L_{\bf m}^{(\nu-1)}(2e^{-i\theta }\cos \theta u)
\Phi _{\bf m}(e^{i\theta }w).$$
Hence the generating formula follows from part (ii) in Proposition \ref{LOG}.
\qed

b) By Lemma \ref{Fourier} we obtain the following evaluation of the spherical Fourier transform of the Laguerre functions
$\Psi _{\bf m}^{(\nu ,\theta )}$: 
$${\cal F}_{\nu } (\Psi _{\bf m}^{(\nu ,\theta )})({\bf s})
=Q_{\bf m}^{(\nu ,\theta )}({\bf s}).$$

We will need a Gutzmer formula for the spherical Fourier transform on a symmetric cone. Let us first state the following Gutzmer formula for the Mellin transform. 

\begin{pro}\label{Gutzmer1}
Let $\psi $ be holomorphic in the following open set in $\mathbb{C}$:
$$\{\zeta =re^{i\theta }\mid r>0,\ |\theta |<\theta _0\}\quad \bigl(0< \theta _0< {\pi \over 2}\bigr).$$
The Mellin transform of $\psi $ is defined by
$${\cal M} \psi (s)=\int _0^{\infty } \psi (r)r^{s-1}dr.$$
Assume that there is a constant $M>0$ such that, for $|\theta |< \theta _0$,
$$\int _0^{\infty } |\psi (re^{i\theta })|^2r^{-1}dr\leq M.$$
Then
$$\int _0^{\infty } |\psi (re^{i\theta })|^2r^{-1}dr
={1\over 2\pi } \int _{\mathbb{R}} |{\cal M}\psi (i\lambda )|^2e^{2\theta \lambda }d\lambda .$$
\end{pro}

Using the decomposition of the symmetric cone $\Omega $ as
$\Omega =]0,\infty [\times \Omega _1$,
where $\Omega _1=\{u\in \Omega \mid \Delta (u)=1\}$, 
one gets the following Gutzmer formula for~$\Omega $:

\begin{pro}\label{Gutzmer2}
Let $\psi $ be a holomorphic function in the tube $T_{\Omega }=\Omega +iV$.
Assume that there are constants $M>0$ and $0<\theta _0< {\pi \over 2}$ such that, for
$|\theta |< \theta _0$,
$$\int _{\Omega }|\psi (e^{i\theta }u)|^2\Delta (u)^{-{N\over n}}m(du)\leq M.$$
Then, for $|\theta |< \theta _0$,
$$\int _{\Omega } |\psi (e^{i\theta }u)|^2\Delta (u)^{-{N\over n}}du 
={1\over (2\pi )^n}\int _{\mathbb{R}^n}|{\cal F}\psi (i\lambda )|^2
e^{2\theta (\lambda _1+\cdots +\lambda _n)}{1\over |c(i\lambda )|^2}m(d\lambda ).$$
\end{pro} 
From Proposition \ref{LOG2} and Proposition \ref{Gutzmer2} we obtain parts (i) and (ii) of Theorem \ref{MPOG2}.
A more general Gutzmer formula has been established for the spherical Fourier transform on Riemannian symmetric spaces of noncompact type \cite{F2004}.
\section{Determinantal formulae}
In the case $d=2$, i.e. $V=Herm(n,\mathbb{C})$, $K=U(n)$, there are determinantal formulae
for the multivariate Laguerre functions
$\Psi _{\bf m}^{(\nu )}$ and for the multivariate Meixner-Pollaczek polynomials $Q_{\bf m}^{(\nu ,\theta )}$.
Consider a Jordan frame $\{c_1,\ldots ,c_n\}$ in $V$, and let $\delta =(n-1,n-2,\ldots ,1,0)$.

\begin{thm} \label{DFL}
Assume $d=2$. 
The multivariate Laguerre function $\Psi _{\bf m}^{(\nu )}$
admits the following determinantal formula involving the one variable Laguerre functions
$\psi _m^{(\nu )}$: 
for $u=\sum _{j=1}^nu_ic_i$,
$$\Psi _{\bf m}^{(\nu )}(u)=\delta ! 2^{-{1\over 2} n(n-1)}
{\det \bigl(\psi _{m_j+\delta _j}^{(\nu -n+1)}(u_i)\bigr)_{1\leq i,j\leq n}
\over V(u_1,\ldots ,u_n)},$$
where $V$ denote the Vandermonde polynomial:
$$V(u_1,\ldots ,u_n)=\prod _{i<j}(u_j-u_i)\ and\ \delta !=\prod _{i=1}^n (n-i)!.$$
As a result one obtains the following determinantal formula for the multivariate
Laguerre polynomials:
$${\bf L}^{\nu }_{\bf m}(u)
=\delta !{\det \Bigl( L_{m_j+\delta _j}^{(\nu -n+1)}(u_i)\Bigr)\over V(u_1,\ldots ,u_n)}.$$
\end{thm}

\proof
We start from the generating formula for the multivariate Laguerre functions (Proposition \ref{LOG}):

\begin{eqnarray*}
{\cal G}_{\nu }^{(3)}(u,w)
&=&\sum _{\bf m} d_{\bf m}\Phi _{\bf m}(w)\Psi _{\bf m}^{(\nu )}(u)\\
&=&\Delta (e-w)^{-\nu }\int _Ke^{-\bigl(ku|(e+w)(e-w)^{-1}\bigr)}dk.
\end{eqnarray*}

In the case $d=2$, the evaluation of this integral is classical:
for $x=\sum _{i=1}^n x_ic_i$, $y=\sum _{j=1}^n y_jc_j$, then
$${\cal I}(x,y)=\int _K e^{(kx|y)}dk
=\delta ! {\det \bigl(e^{x_iy_j}\bigr)\over V(x_1,\ldots ,x_n)V(y_1,\ldots ,y_n)}.$$
Therefore, for $u=\sum _{i=1}^n u_ic_i$, $w=\sum _{j=1}^n w_jc_j$,
$${\cal G}_{\nu }^{(3)}(u,w)=
\delta ! \prod _{j=1}^n (1-w_j)^{-\nu }
{\det \Bigl(e^{-u_i{1+w_j\over 1-w_j}}\Bigr)
\over V(u_1,\ldots ,u_n)V\bigl({1+w_1\over 1-w_1},\ldots ,{1+w_n\over 1-w_n}\bigr)}.$$
Noticing that
$${1+w_j\over 1-w_j}-{1+w_k\over 1-w_k}=2{w_j-w_k\over (1+w_j)(1+w_k)},$$
we obtain
$${\cal G}_{\nu }^{(3)}(u,w)=\delta ! 2^{-{1\over 2} n(n-1)}
{\det \Bigl((1-w_j)^{-(\nu-n+1)}e^{-u_i{1+w_j\over 1-w_j}}\Bigr)\over
V(u_1,\ldots , u_n)V(w_1,\ldots ,w_n)}.$$
We will expand the above expression in Schur function series by using a formula due to Hua 
(See \cite{H1963}, Theorem 1.2.1, p.22).

\begin{lem}\label{HuaTheorem}
Consider $n$ power series 
$$f_i(w)=\sum _{m=0}^{\infty } c_m^{(i)}w^m\quad (i=1,\ldots ,n).$$
Then
$${\det \bigl( f_i(w_j)\bigr)\over V(w_1,\ldots ,w_n)}
=\sum _{\bf m}a_{\bf m}s_{\bf m}(w_1,\ldots ,w_n) ,$$
where $s_{\bf m}$ is the Schur function associated to the partition $\bf m$, and
$$a_{\bf m}=\det \bigl(c_{m_j+\delta _j}^{(i)}\bigr).$$
\end{lem}


Let $\nu '=\nu -n+1$, and consider the $n$  power series
$$f_i(w):=(1-w)^{-\nu '}e^{-u_i{1+w\over 1-w}}
=\sum _{m=0}^{\infty } \psi _m^{(\nu ')}(u_i)w^m.$$
Since
$$d_{\bf m}\Phi _{\bf m}\Bigl(\sum _{j=1}^n w_jc_j\Bigr)=s_{\bf m}(w_1,\ldots ,w_n),$$
we obtain
$$\Psi _{\bf m}^{(\nu )}(u)=\delta ! 2^{-{1\over 2} n(n-1)}
{\det \bigl(\psi _{m_j+\delta _j}^{(\nu -n+1)}(u_i)\bigr)\over V(u_1,\ldots ,u_n)}.$$
\qed

By using the same method we will obtain a determinantal formula for the multivariate 
Meixner-Pollaczek polynomials $Q_{\bf m}^{(\nu ,\theta )}$.

\begin{thm}\label{DFMP}
Assume $d=2$. Then
$$Q_{\bf m}^{(\nu ,\theta )}({\bf s})
=(-2\cos \theta )^{-{1\over 2}n(n-1)} \delta !
{\det \Bigl(q_{m_j+\delta _j}^{(\nu-n+1,\theta )}(s_i)\Bigr)_{1\leq i,j\leq n}
\over V(s_1,\ldots ,s_n)},$$
where $q_m^{(\nu ,\theta )}$ denotes the one variable Meixner-Pollaczek polynomial.
\end{thm}

\proof
We start from the generating formula for the multivariate Meixner-Pollaczek polynomials
$Q_{\bf m}^{(\nu ,\theta )}$
(Theorem \ref{MPOG2}, (ii)):
$$\sum _{\bf m} d_{\bf m}Q_{\bf m}^{(\nu ,\theta )}({\bf s})\Phi _{\bf m}(w)
=\Delta \bigl((e-e^{i\theta }w)(e+e^{-i\theta }w )\bigr)^{-{\nu \over 2}}
\varphi _{\bf s}\bigl(c_{\theta }(w)^{-1}\bigr) .$$
For $x=\sum _{i=1}^n x_ic_i$, the spherical function $\varphi _{\bf s}(x)$ is essentially
a Schur function in the variables $x_1,\ldots ,x_n$:
$$\varphi _{\bf s}(x)=\delta !(x_1x_2\ldots x_r)^{{1\over 2}(n-1)}
{\det (x_j^{s_i})\over V(s_1,\ldots ,s_n)V(x_1,\ldots ,x_n)}.$$
Let us compute now, for $w=\sum _{j=1}^n w_jc_j$,
\begin{eqnarray*}
& & \Delta \bigl((e-e^{i\theta }w)(e+e^{-i\theta }w )\bigr)^{-{\nu \over 2}}
\varphi _{\bf s}\bigl(c_{\theta }(w)^{-1}\bigr) \\
& &=\delta ! \prod _{j=1}^n (1-2i\sin \theta w_j-w_j^2)^{-{\nu \over 2} } \\
& &\times \prod _{j=1}^n\Bigl(c_{\theta }(w_j)\Bigr)^{{1\over 2}(n-1)} 
{\det \Bigl(\bigl(c_{\theta }(w_j)\bigr)^{-s_i}\Bigr)\over 
V(s_1,\ldots ,s_n)V\Bigl(c_{\theta }(w_1),\ldots ,c_{\theta }(w_n)\Bigr)}.
\end{eqnarray*}
In the same way, as for the proof of Theorem \ref{DFL},
we obtain
\begin{eqnarray*}
& &\Delta \bigl((e-e^{i\theta }w)(e+e^{-i\theta })\bigr)^{-{\nu \over 2}}
\varphi _{\bf s}\bigl(c_{\theta }(w)^{-1}\bigr) \\ 
&=&(-2\cos \theta )^{-{1\over 2}n(n-1)}\delta ! 
{\det \Bigl((1-e^{i\theta }w_j)^{s_i-{\nu \over 2}+{1\over 2}(n-1)}
(1+e^{-i\theta }w_j)^{-s_i-{\nu \over 2}+{1\over 2}(n-1)}\Bigr)
\over V(s_1,\ldots ,s_n)V(w_1,\ldots ,w_n)}.
\end{eqnarray*}
We apply once more Lemma \ref{HuaTheorem} to the $n$ power series
$$f_i(w):=(1-e^{i\theta }w)^{s_i-{\nu '\over 2}}(1+e^{-i\theta}w)^{-s_i-{\nu '\over 2}}
=\sum _m^{\infty }q_m^{(\nu ', \theta )}(s_i)w^m$$
with $\nu '=\nu -n+1$,
and obtain finally:
$$Q_{\bf m}^{(\nu ,\theta )}({\bf s})
=(-2\cos \theta )^{-{1\over 2}n(n-1)}\delta !
{\det \Bigl(q_{m_j+\delta _j}^{(\nu-n+1,\theta )}(s_i)\Bigr)
\over V(s_1,\ldots ,s_n)}.$$
\qed

\section{Difference equation for the Meixner-Pollaczek polynomials $Q_{\bf m}^{(\nu ,\theta )}$}
The one variable Meixner-Pollaczek polynomials $q_m=q_m^{(\nu ,\theta )}$ satisfies the following difference
equation
$$e^{-i\theta }\Bigl(s+{\nu \over 2}\Bigr)\bigl(q_m(s+1)-q_m(s)\bigr)
+e^{i\theta }\Bigl(-s+{\nu \over 2}\Bigr)\bigl(q_m(s-1)-q_m(s)\bigr)
=2m\cos \theta q_m.$$
(See \cite{AAR1999}, p.348, 37.(d)).
We will establish an analogue of this formula for the multivariate Meixner-Pollaczek polynomials
$Q_{\bf m}^{(\nu ,\theta )}$.

Recall the Pieri's formula for the spherical functions:
$${\rm tr}\,  u\, \varphi _{\bf s}(u)=\sum _{j=1}^n \alpha _j({\bf s})\varphi _{{\bf s}+\varepsilon _j}(u),\
{\rm with}\  \alpha _j({\bf s})=\prod _{k\not= j}{s_j-s_k+{d\over 2}\over s_j-s_k}$$
($\{\varepsilon _i\}$ denotes the canonical basis of $\mathbb{C}^n$). 
See \cite{Di1990}, Proposition 6.1 or \cite{Z1995}, Theorem 1, and also \cite{L1998}, p.320.
We introduce the difference operator $D_{\nu ,\theta }$:
\begin{eqnarray*}
D_{\nu ,\theta }f({\bf s})
& &=e^{-i\theta }\sum _{j=1}^n\bigl(s_j+{\nu \over 2}-{d\over 4}(n-1)\bigr) \alpha _j({\bf s})
\bigl(f({\bf s}+\varepsilon _j)-f({\bf s})\bigr) \\
& &+e^{i\theta }\sum _{j=1}^n\bigl(-s_j+{\nu \over 2}-{d\over 4}(n-1)\bigr)
\alpha _j(-{\bf s})\bigl(f({\bf s}-\varepsilon _j)-f({\bf s})\bigr).
\end{eqnarray*}

\begin{thm}\label{DifferenceMP}
The Meixner-Pollaczek polynomial $Q_{\bf m}^{(\nu ,\theta )}$ is an eigenfunction of the difference operator $D_{\nu ,\theta }$:
$$D_{\nu ,\theta }Q_{\bf m}^{(\nu ,\theta )}=2|{\bf m}|\cos \theta \ Q_{\bf m}^{(\nu ,\theta )}.$$
\end{thm}

\bigskip

For the proof we will use the scheme we have used in the proof of Theorem \ref{MPOG1}.
For $i=1,2,3,4$, we define the operators $D_{\nu ,\theta }^{(i)}$.
The operator $D_{\nu ,\theta }^{(1)}=D_{\theta }^{(1)}$ is a first order differential  operator on the domain $\cal D$:
$$D_{\theta }^{(1)}f=e^{i\theta }\langle w+e,\nabla f\rangle
+e^{-i\theta }\langle w-e,\nabla f\rangle .$$
(For $w_1,w_2\in V_\mathbb{C}$, $\langle w_1,w_2\rangle ={\rm tr}\,  (w_1w_2)$.)
The operators $D_{\nu ,\theta }^{(i)}$, for $i=2,3,4$ are defined by the relations:
$$D_{\nu ,\theta }^{(2)}C_{\nu } =C_{\nu }D_{\nu ,\theta }^{(1)},\
{\cal L}_{\nu } D_{\nu ,\theta }^{(3)}=D_{\nu ,\theta }^{(2)}{\cal L}_{\nu },\
{\cal F}_{\nu } D_{\nu ,\theta }^{(3)}=D_{\nu ,\theta }^{(4)}{\cal F}_{\nu }.$$
The operator $D_{\nu ,\theta }^{(2)}$ is a first order differential operator on the tube $T_{\Omega }$.
In Section~8 we will see that $D_{\nu ,\theta }^{(3)}$ is a second order differential operator on the cone $\Omega $,
and prove that $D_{\nu ,\theta }^{(4)}$ is the difference operator $D_{\nu ,\theta }$ we have introduced
above.

The function
$\Phi _{\bf m}^{(\theta )}(w)=\Phi _{\bf m}(w\cos \theta +ie\sin \theta )$
is an eigenfunction of the operator $D_{\theta }^{(1)}$:
$D_{\theta }^{(1)}\Phi _{\bf m}^{(\theta )}=2|{\bf m}|\cos \theta \ \Phi _{\bf m}^{(\theta )}$.
Hence $F_{\bf m}^{(\nu ,\theta )}=C_{\nu } \Phi _{\bf m}^{(\theta )}$ is an eigenfunction of
$D_{\nu ,\theta }^{(2)}$:
$D_{\nu ,\theta }^{(2)}F_{\bf m}^{(\nu ,\theta )}=2|{\bf m}|\cos \theta \ F_{\bf m}^{(\nu ,\theta )}$.
Further, since 
${\cal L}_{\nu } \Psi _{\bf m}^{(\nu ,\theta )}=
{(\nu )_{\bf m}\over \bigl({N\over n}\bigr)_{\bf m}}F_{\bf m}^{(\nu ,\theta )}$, we get
$D_{\nu ,\theta }^{(3)}\Psi _{\bf m}^{(\nu ,\theta )}=2|{\bf m}|\cos \theta \ \Psi _{\bf m}^{(\nu ,\theta )}$.
Finally, since $Q_{\bf m}^{(\nu ,\theta )}={\cal F}_{\nu } \Psi _{\bf m}^{(\nu ,\theta )}$, then
$D_{\nu ,\theta }^{(4)}Q_{\bf m}^{(\nu ,\theta )}=2|{\bf m}|\cos \theta \ Q_{\bf m}^{(\nu ,\theta )}$.
Hence the proof of Theorem \ref{DifferenceMP} amounts to showing that $D_{\nu ,\theta }^{(4)}=D_{\nu ,\theta }$.
\section{The symmetries $S_{\nu }^{(i)}$ ($i=1,2,3,4$) and the Hankel transform}
The symmetries  $S_{\nu }^{(i)}$ we introduce now will be useful for the computation of the operators $D_{\nu ,\theta }^{(i)}$.
We start from the symmetry $w\mapsto -w$ of the domain $\cal D$. Its action on functions is given by
$S^{(1)}f(w)=f(-w)$.
We carry this symmetry over the tube $T_{\Omega }$ through the Cayley transform and obtain the inversion $z\mapsto z^{-1}$. We define $S_{\nu }^{(2)}$ such that
$S_{\nu }^{(2)}C_{\nu } =C_{\nu } S^{(1)}$.
Hence, for a function $F$ on $T_{\Omega }$,
$S_{\nu }^{(2)}F(z)=\Delta (z)^{-\nu }F(z^{-1})$.
Further $S_{\nu }^{(3)}$ is defined by the relation
${\cal L}_{\nu } S_{\nu }^{(3)}=S_{\nu }^{(2)}{\cal L}_{\nu }$.
By a generalized Tricomi theorem (Theorem XV.4.1 in \cite{FK1994}), the unitary isomorphism 
$S_{\nu }^{(3)}$ of $L_{\nu }^2(\Omega )$ is the Hankel transform: $S_{\nu }^{(3)}=U_{\nu }$,
$$U_{\nu } \psi (u)=\int _{\Omega } H_{\nu } (u,v)\psi (v)\Delta (v)^{\nu -{N\over n}}m(dv).$$
The kernel $H_{\nu }(u,v)$ has the following invariance property:  for $g\in G$,
$$H_{\nu } (g\cdot u,v)=H_{\nu } (u, g^*\cdot v),\ 
{\rm and}\ 
H_{\nu }(u,e)={1\over \Gamma _{\Omega }(\nu )}{\cal J}_{\nu } (u),$$
where ${\cal J}_{\nu }$ is a multivariate Bessel function.

Finally we define $S_{\nu }^{(4)}$ acting on symmetric polynomials in $n$ variables such that
$$S_{\nu }^{(4)} {\cal F}_{\nu } ={\cal F}_{\nu } S_{\nu } ^{(3)}.$$

\begin{pro}\label{Symmetry}
For a function $\psi $ on $\Omega $ of the form
$\psi (u)=e^{-{\rm tr} u}q(u)$,
where $q$ is a $K$-invariant polynomial,
${\cal F}_{\nu }(U_{\nu }\psi )({\bf s})={\cal F}_{\nu }\psi (-{\bf s})$.
It follows that, for a symmetric polynomial $p$ on $\mathbb{C}^n$,
$$S_{\nu }^{(4)}p({\bf s})=p(-{\bf s}).$$
\end{pro}

\proof
We will evaluate the spherical Fourier transform ${\cal F}_{\nu } (U_{\nu }\psi )$.
By the invariance property, the kernel $H_{\nu }(u,v)$ can be written
$$H_{\nu }(u,v)=h_{\nu } \bigl(P(v^{1\over 2})u\bigr)\Delta (u)^{-{\nu \over 2}}\Delta (v)^{-{\nu \over 2}},$$
with $h_{\nu }(u)=H_{\nu }(u,e)\Delta (u)^{\nu \over 2}$,
and $P$ is the so-called quadratic representation of the Jordan algebra $V$.
Let us compute first
\begin{eqnarray*}
& &\int _{\Omega } H_{\nu }(u,v)\varphi _{\bf s}(u)\Delta (u)^{{\nu \over 2}-{N\over n}}m(du) \\
&=& \Delta (v)^{-{\nu \over 2}}\int _{\Omega }h_{\nu }\bigl(P(v^{1\over 2})u\bigr)\varphi _{\bf s}(u)\Delta (u)^{-{N\over n}}m(du).
\end{eqnarray*}
By letting $P(v^{1\over 2})u=u'$, we get
\begin{eqnarray*}
& &\int _{\Omega } H_{\nu }(u,v)\varphi _{\bf s}(u)\Delta (u)^{{\nu \over 2}-{N\over n}}m(du) \\
&=&\Delta (v)^{-{\nu \over 2}}\int _{\Omega } h_{\nu } (u')\varphi _{\bf s}\bigl(P(v^{-{1\over 2}})u'\bigr)
\Delta (u')^{-{N\over n}}m(du').
\end{eqnarray*}
By using $K$-invariance and the functional equation of the spherical function $\varphi _{\bf s}$,
$$\int _K\varphi _{\bf s}\bigl(P(v^{-{1\over 2}})ku')dk=\varphi _{\bf s}(v^{-1})\varphi _{\bf s}(u'),$$
we get
$$\int _{\Omega } H_{\nu }(u,v)\varphi _{\bf s}(u)\Delta (u)^{{\nu \over 2}-{N\over n}}m(du)
=\varphi _{\bf s}(v^{-1})\Delta (v)^{-{\nu \over 2}}{\cal F}(h_{\nu } )({\bf s}).$$
Recall that $\varphi _{\bf s}(v^{-1})=\varphi _{-{\bf s}}(v)$.
We multiply both sides by $\psi (v)$ and get by integrating with respect to $v$:
$$\Gamma _{\Omega }\bigl({\bf s}+{\nu \over 2}+\rho \bigr)
{\cal F}_{\nu } (U_{\nu } \psi )({\bf s})
={\cal F}h_{\nu }({\bf s})\Gamma _{\Omega }\bigl(-{\bf s}+{\nu \over 2}+\rho \bigr)
{\cal F}_{\nu } \psi (-{\bf s}).$$
Consider the special case $\psi (u)=\Psi _0(u)=e^{-{\rm tr}\,  u}$. Since
$U_{\nu }\Psi _0=\Psi _0$, and ${\cal F}_{\nu } \Psi _0 \equiv 1$, we get
$${\cal F}(h_{\nu })({\bf s})={\Gamma _{\Omega }\bigl({\bf s}+{\nu \over 2}+\rho \bigr)
\over \Gamma _{\Omega }\bigl(-{\bf s}+{\nu \over 2}+\rho \bigr)}.$$
Finally
${\cal F}_{\nu } (U_{\nu }\psi )({\bf s})={\cal F}_{\nu } \psi (-{\bf s}),$
and $S_{\nu } ^{(4)}p({\bf s})=p(-{\bf s})$.
\qed

\begin{cor}\label{Q-symmetry}
$$Q_{\bf m}^{(\nu ,\theta )}(-{\bf s})=(-1)^{|{\bf m}|}Q_{\bf m}^{(\nu,-\theta )}({\bf s}).$$
\end{cor}

\proof
This relation follows from
$$S^{(1)}\Phi _{\bf m}^{(\theta )}=\Phi _{\bf m}^{(\theta )}(-w)
=(-1)^{|{\bf m}|}\Phi _{\bf m}^{(-\theta )}(w),$$
which is easy to check, and Proposition \ref{Symmetry}.
\qed

The operator $D_{\nu ,\theta }^{(i)}$ ($i=1,2,3,4$) can be written
$$D_{\nu ,\theta }^{(i)}=e^{i\theta }D_{\nu }^{(i,+)}+e^{-i\theta }D_{\nu }^{(i,-)}.$$
For $i=1$, $D_{\nu }^{(1,\pm )}$ does not depend on $\nu $, $D_{\nu }^{(1,\pm )}=D^{(1,\pm)}$:
$$D^{(1,+)}f(w)=\langle w+e,\nabla f (w)\rangle,\quad
D^{(1,-)}f(w)=\langle w-e,\nabla f(w)\rangle .$$
Observe that
$D^{(1,-)}=S^{(1)}D^{(1,+)}S^{(1)}$.
Hence, for $i=2,3,4$,
$D_{\nu }^{(i,-)}=S_{\nu }^{(i)}D_{\nu }^{(i,+)}S_{\nu }^{(i)}$.
In next Section we will compute first $D_{\nu }^{(i,-)}$. The operator $D_{\nu }^{(i,+)}$
is then obtained by using the above relation. 
For $i=3$, we will use the following property of the Hankel transform

\begin{pro}
$$U_{\nu } ({\rm tr}\,  v\ \psi )
=-\Bigl(\langle u,\bigl({\partial \over \partial u}\bigr)^2\rangle 
+\nu {\rm tr}\,  \bigl({\partial \over \partial u}\bigr)\Bigr)U_{\nu } \psi .$$
\end{pro}

This is a consequence of Proposition XV.2.3 in \cite{FK1994}.
 \section{Proof of Theorem \ref{DifferenceMP}}
a) Recall that $D^{(1,-)}$ is the first order differential operator on the domain $\cal D$ given by
$$D^{(1,-)}f(w)=\langle w-e,\nabla f(w)\rangle ,$$
and $D_{\nu }^{(2,-)}$ is the first order differential operator on the tube $T_{\Omega }$ such that
$$D_{\nu } ^{(2,-)}C_{\nu } =C_{\nu } D^{(1,-)}.$$

\begin{lem}\label{DiffrenceF}
$$D_{\nu } ^{(2,-)}F(z)=-\langle z+e,\nabla F(z)\rangle -n\nu F(z).$$
\end{lem}

\proof
Recall that, for a function $F$ on the tube $T_{\Omega }$,
$$f(w)=(C_{\nu }^{-1}F)(w)=\Delta (e-w)^{-\nu } F\bigl(c(w)\bigr),$$
where $c$ is the Cayley transform
$$c(w)=(e+w)(e-w)^{-1}=2(e-w)^{-1}-e.$$
Its differential is given by
$$(Dc)_w=2P\bigl((e-w)^{-1}\bigr).$$
We get
$$\nabla f(w)=\nabla \bigl(\Delta (e-w)^{-\nu }\bigr)F\bigl(c(w)\bigr)
+\Delta (e-w)^{-\nu } 2P\bigl(e-w)^{-1}\bigr)\Bigl(\nabla F\bigl(c(w)\bigr)\Bigr).$$
By using $\nabla \bigl(\Delta (x)^{\alpha } \bigr)=\alpha \Delta (x)^{\alpha } x^{-1},$ and
$$\langle e-w,(e-w)^{-1}\rangle =n,\ 
P\bigl((e-w)^{-1}\bigr)(e-w)=(e-w)^{-1},$$
we obtain
\begin{eqnarray*}
& &D^{(1,-)}f(w)=\langle w-e,\nabla f(w)\rangle \\
&= &\Delta (e-w)^{-\nu } \Bigl(-n\nu F\bigl(c(w)\bigr)
+2\langle (w-e)^{-1},\nabla F\bigl(c(w)\bigr)\rangle \Bigr)
=(C_{\nu }^{-1}G)(z),
\end{eqnarray*}
with
$$G(z)=-\langle z+e,\nabla F(z)\rangle -n\nu F(z).$$
\qed

b) Consider now the differential operator $D_{\nu }^{(3,-)}$ on the cone $\Omega $ such that
$${\cal L}_{\nu } D_{\nu }^{(3,-)}=D_{\nu }^{(2,-)}{\cal L}_{\nu }.$$
Recall that  the modified Laplace transform ${\cal L}_{\nu } \psi $ of a function $\psi $, defined on $\Omega $, is given by
$$F(z)={\cal L}_{\nu } \psi (z)
={2^{n\nu } \over \Gamma _{\Omega } (\nu )}
\int _{\Omega } e^{-(z|u)}\psi (u)\Delta (u)^{\nu -{N\over n}}m(du).$$

\begin{lem}\label{DifferencePhi}
$$D_{\nu } ^{(3,-)}\psi (u)=\langle u,\nabla \psi (u)\rangle +{\rm tr}\,  u\, \psi (u).$$
\end{lem}

\proof
For $a\in V_\mathbb{C}$,
$$\langle a,\nabla F(z)\rangle 
={2^{n\nu }\over \Gamma _{\Omega } (\nu )}\int _{\Omega }
e^{-(z|u)}(-\langle a,u\rangle )\psi (u)\Delta (u)^{\nu -{N\over n}}m(du).$$
Observe that $(z|u)e^{-(z|u)}=\langle u,\nabla _u\rangle e^{-(z|u)}$. 
Therefore
$$\langle z,\nabla F(z)\rangle
={2^{n\nu }\over \Gamma _{\Omega }(\nu )}\int _{\Omega }
(-\langle u,\nabla _u\rangle e^{-(z|u)})\psi (u)\Delta (u)^{\nu -{N\over n}}m(du).$$
An integration by parts gives
$$={2^{n\nu }\over \Gamma _{\Omega }(\nu )}\int _{\Omega }
e^{-(z|u)}(\langle u,\nabla \rangle +n\nu )\psi (u)\Delta ^{\nu -{N\over n}}m(du).$$
Finally
$$(D_{\nu }^{(2,-)}F)(z)={\cal L}_{\nu } (\langle u,\nabla \psi \rangle +{\rm tr}\,  u \, \psi ).$$
\qed

c) The operator $D_{\nu }^{(4,-)}$ acting on symmetric functions on $\mathbb{C}^n$ is such that
$$D_{\nu }^{(4,-)}{\cal F}_{\nu } ={\cal F}_{\nu } D_{\nu }^{(3,-)}.$$
Recall that the spherical Fourier transform $f={\cal F}_{\nu } \psi $  of a function $\psi $, defined  on $\Omega $,
is given by
$$f({\bf s})=({\cal F}_{\nu } \psi )({\bf s})
={1\over \Gamma _{\Omega } \bigl({\bf s}+{\nu \over 2}+\rho \bigr)}
\int _{\Omega } \varphi _{\bf s}(u)\psi (u)\Delta (u)^{{\nu \over 2}-{N\over n}}m(du).$$  

\begin{pro}\label{Difference4}
The operator $D_{\nu } ^{(4,-)}$ is the following difference operator: 
for a function $f$ on $\mathbb{C}^n$,
$$D_{\nu } ^{(4,-)}f({\bf s})
=\sum _{j=1}^n \bigl(s_j+{\nu \over 2}-{d\over 4}(n-1)\alpha _j({\bf s})\bigr)
\bigl(f({\bf s}+\varepsilon _j)-f({\bf s})\bigr).$$
\end{pro}

\proof
We will compute 
${\cal F}_{\nu }(D_{\nu }^{(3,-)}\psi )={\cal F}_{\nu }(\langle u,\nabla \psi \rangle +{\rm tr}\,  u\, \psi )$.
Consider first
$${\cal F}_{\nu } (\langle u,\nabla \psi \rangle )({\bf s})
={1\over \Gamma _{\Omega } \bigl({\bf s}+{\nu \over 2}+\rho \bigr)}
\int _{\Omega }\langle u,\nabla \psi (u)\rangle \varphi _{{\bf s}+{\nu \over 2}}(u)
\Delta (u)^{-{N\over n}}m(du).$$
An integration by parts gives,
since the function $\varphi _{\bf s}$ is homogeneous of degree $\sum _{j=1}^n s_j$
(observe that $\sum _{j=1}^n\rho _j=0$),
\begin{eqnarray*}
&=&{1\over \Gamma _{\Omega } \bigl({\bf s}+{\nu \over 2}+\rho \bigr)}
\int _{\Omega } \psi (u)\bigl(-\langle u,\nabla _u\rangle 
\varphi _{{\bf s}+{\nu \over 2}}(u)\bigr)\Delta (u)^{-{N\over n}}m(du)\\
&=&{1\over \Gamma _{\Omega } \bigl({\bf s}+{\nu \over 2}+\rho \bigr)}
\int _{\Omega } \psi (u)\Bigl(-\sum _{j=1}^n \bigl(s_j+{\nu \over 2}\bigr)\Bigr)
\varphi _{\bf s}(u)\Delta (u)^{{\nu \over 2}-{N\over n}}m(du)\\
&=&-\sum _{j=1}^n \bigl(s_j+{\nu \over 2}\bigr){\cal F}_{\nu } \psi ({\bf s}).
\end{eqnarray*}
Recall the Pieri's formula for the spherical functions:
$${\rm tr}\,  u\, \varphi _{\bf s}(u)=\sum _{j=1}^n \alpha _j({\bf s})\varphi _{{\bf s}+\varepsilon _j}(u),\
{\rm with}\  \alpha _j({\bf s})=\prod _{k\not= j}{s_j-s_k+{d\over 2}\over s_j-s_k}.$$
Hence
\begin{eqnarray*}
& &{\cal F}_{\nu } ({\rm tr}\,  u\, \psi )({\bf s}) \\
&=&{1\over \Gamma _{\Omega } \bigl({\bf s}+{\nu \over 2}+\rho \bigr)}
\int _{\Omega }\psi (u)\Bigl(\sum _{j=1}^n \alpha ({\bf s})\varphi _{{\bf s}+\varepsilon _j}(u)\Bigr)
\Delta (u)^{{\nu \over 2}-{N\over n}}m(du) \\
&=&\sum _{j=1}^n{\Gamma _{\Omega }\bigl({\bf s}+\varepsilon _j+{\nu \over 2}+\rho \bigr)
\over  \Gamma _{\Omega } \bigl({\bf s}+{\nu \over 2}+\rho \bigr)} \alpha _j({\bf s}) \\
& &\times {1\over \Gamma _{\Omega }\bigl({\bf s}+\varepsilon _j+{\nu \over 2}+\rho \bigr)}
\int _{\Omega }\psi (u)\varphi _{{\bf s}+\varepsilon _j}(u)\Delta ^{{\nu \over 2}-{N\over n}}m(du) \\
&=&\sum _{j=1}^n\bigl(s_j+{\nu \over 2}-{d\over 4}(n-1)\bigr)\alpha _j({\bf s})
{\cal F}_{\nu } \psi ({\bf s}+\varepsilon _j).
\end{eqnarray*}
Finally
\begin{eqnarray*}
& &{\cal F}_{\nu } (D_{\nu }^{(3,-)}\psi )({\bf s}) \\
&=&\sum _{j=1}^n \bigl(s_j+{\nu \over 2}-{d\over 4}(n-1)\bigr)\alpha _j({\bf s})f({\bf s}+\varepsilon _j)
-\sum _{j=1}^n \bigl(s_j+{\nu \over 2}\bigr)f({\bf s}),
\end{eqnarray*}
with $f={\cal F}_{\nu } (\psi )$. From $D_{\nu }^{(3,-)}\Psi _0=0$ and ${\cal F}_{\nu } (\Psi _0)=1$, we get
$$\sum _{j=1}^n \bigl(s_j+{\nu \over 2}-{d\over 4}(n-1)\bigr)\alpha _j({\bf s})
=\sum _{j=1}^n\bigl(s_j+{\nu \over 2}\bigr).$$
Therefore
$${\cal F}_{\nu } (D_{\nu }^{(3,-)}\psi )({\bf s}) 
=\sum _{j=1}^n \bigl(s_j+{\nu \over 2}-{d\over 4}(n-1)\bigr)
\alpha _j({\bf s})\bigl(f({\bf s}+\varepsilon _j)-f({\bf s})\bigr).$$
\qed

We finish now the proof of Theorem \ref{DifferenceMP}. Recall that
$$D_{\nu }^{(4,+)}=S_{\nu }^{(4)}D_{\nu }^{(4,-)}S_{\nu }^{(4)},\quad 
{\rm and} \quad S_{\nu }^{(4)}f({\bf s})=f(-{\bf s}).$$
Therefore, by Proposition \ref{Difference4},
$$D_{\nu }^{(4,+)}f({\bf s}) 
=\sum _{j=1}^n \bigl(-s_j+{\nu \over 2}-{d\over 4}(n-1)\bigr)
\alpha _j(-{\bf s})\bigl(f({\bf s}-\varepsilon _j)-f({\bf s})\bigr).$$
We have established the formula of Theorem \ref{DifferenceMP} since
$$D_{\nu ,\theta }=D_{\nu ,\theta }^{(4)}=e^{i\theta }D_{\nu }^{(4,+)}+e^{-i\theta }D_{\nu }^{(4,-)}.$$
\section{Pieri's formula for the Meixner-Pollaczek polynomials $Q_{\bf m}^{(\nu ,\theta )}$}
\begin{thm}\label{Pieri}
The Meixner-Pollaczek polynomials $Q_{\bf m}^{(\nu ,\theta )}$ satisfy the following Pieri's formula:
\begin{eqnarray*}
& &(2|{\bf s}|\cos \theta -2i|2{\bf m}+\nu |\sin \theta )Q_{\bf m}^{(\nu ,\theta )}({\bf s}) \\
&=&\sum _{j=1}^n\bigl(m_j+\nu -1-{d\over 4}(j-1)\bigr)\alpha _j({\bf m}-\varepsilon _j-\rho )
d_{{\bf m}-\varepsilon _j}Q_{{\bf m}-\varepsilon _j}^{(\nu ,\theta )}({\bf s}) \\
& &-\sum _{j=1}^n \bigl(m_j+1+{d\over 4}(n-j)\bigr)\alpha _j(-{\bf m}-\varepsilon _j-\rho )
d_{{\bf m}+\varepsilon _j}Q_{\bf m}^{(\nu ,\theta )}({\bf s}).
\end{eqnarray*}
\end{thm}

\proof
The generating formula (Theorem \ref{MPOG1} (ii)), with ${\bf s}={\bf m}+{\nu \over 2}-\rho $ can be written:
\begin{eqnarray*}
& &\sum _{\bf k} d_{\bf k} Q_{\bf k}^{(\nu ,\theta )}\bigl({\bf m}+{\nu \over 2}-\rho \bigr) \Phi _{\bf k}(w)\\
&=&\Delta (e+e^{-i\theta }w)^{-\nu }\Phi _{\bf m}\bigl((e-e^{i\theta }w)(e+e^{-i\theta } w)^{-1}\bigr).
\end{eqnarray*}
Since
\begin{eqnarray*}
& &F_{\bf m}^{(\nu ,\theta )}(e^{-i\theta }w) \\
&=&2^{n\nu } \Delta (e+e^{-i\theta }w)^{-\nu }
(-1)^{|{\bf m}|} e^{-i|{\bf m}|\theta } \Phi _{\bf m} \bigl((e-e^{i\theta }w)(e+e^{-i\theta }w)^{-1}\bigr),
\end{eqnarray*}
we obtain
$$\sum _{\bf k} Q_{\bf k}^{(\nu ,\theta )}\bigl({\bf m}+{\nu \over 2}-\rho \bigr)
e^{i|{\bf k}|\theta }\Phi _{\bf k}(w)
=2^{-n\nu } (-1)^{|{\bf m}|} e^{i|{\bf m}|\theta }F_{\bf m}^{(\nu ,\theta )}(w).$$
Recall that the function $F_{\bf m}^{(\nu ,\theta )}$ is an eigenfunction of the differential operator 
$D_{\nu ,\theta }^{(2)}$:
$$D_{\nu ,\theta }^{(2)}F_{\bf m}^{(\nu ,\theta )}(w)
=2|{\bf m}|\cos \theta F_{\bf m}^{(\nu ,\theta )}(w).$$
It follows that
\begin{eqnarray*}
& &\sum _{\bf k}d_{\bf k}Q_{\bf k}^{(\nu ,\theta )}\bigl({\bf m}+{\nu \over 2}-\rho \bigr)
e^{i|{\bf k}|\theta }D_{\nu ,\theta }^{(2)}\Phi _{\bf k}(w) \\
&=&2|{\bf m}|\cos \theta \sum _{\bf k}d_{\bf k}Q_{\bf k}^{(\nu ,\theta )}\bigl({\bf m}+{\nu \over 2}-\rho \bigr)
\Phi _{\bf k}(w). \qquad \label ((9.1)
\end{eqnarray*}
\qed

To prove Theorem \ref{Pieri} we will compute $D_{\nu,\theta }^{(2)}\Phi _{\bf k}(w)$.

\begin{lem}\label{TwoFormulas}
The following formulas hold. 
{\rm (i)}
$${\rm tr}\,  \bigl(\nabla \varphi _{\bf s}(z)\bigr)
=\sum _{j=1}^n \bigl(s_j+{d\over 4}(n-1)\bigr)\alpha _j(-{\bf s})\varphi _{{\bf s}-\varepsilon _j}(z).$$
{\rm (ii)}
\begin{eqnarray*}
& &D_{\nu ,\theta }^{(2)}\varphi _{\bf s}(z) 
=e^{i\theta }\Bigl(\sum _{j=1}^n \bigl(s_j-{d\over 4}(n-1)+\nu \bigr)
\alpha _j({\bf s})\varphi _{{\bf s}+\varepsilon _j}(z)+\bigl(\sum _{j=1}^n s_j\bigr)\varphi _{\bf s}(z)\Bigr)\\
& &-e^{-i\theta }\Bigl(\sum _{j=1}^n \bigl(s_j+{d\over 4}(n-1)\bigr)
\alpha _j(-{\bf s})\varphi _{{\bf s}-\varepsilon _j}(z)
+\bigl(\sum _{j=1}^n s_j\bigr)\varphi _{\bf s}(z)+n\nu \varphi _{\bf s}(z)\Bigr).
\end{eqnarray*}
\end{lem}

\proof
(i) For $t>0$ we consider the following Laplace integral:
$$\int _{\Omega }e^{-(x|y)}e^{-t {\rm tr}\, y}\varphi _{\bf s}(y)\Delta (y)^{-{N\over n}}m(dy)
=\Gamma _{\Omega } ({\bf s}+\rho )\varphi _{-{\bf s}}(te+x).$$
Taking the derivatives with respect to $t$ for $t=0$, one gets:
$$-\int _{\Omega }e^{-(x|y)}{\rm tr}\,  y\, \varphi _{\bf s}(y)\Delta (y)^{-{N\over n}}m(dy)
=\Gamma _{\Omega }({\bf s}+\rho )\, \rm{tr}\,  \bigl(\nabla \varphi _{-{\bf s}}(x)\bigr).$$
By using Pieri's formula for the spherical functions,
$${\rm tr}\,  y\, \varphi _{\bf s}(y)=\sum _{j=1}^n\alpha _j({\bf s})\varphi _{{\bf s}+\varepsilon _j}(y),$$
and since
$$\sum _{j=1}^n \alpha _j({\bf s})\int _{\Omega } 
e^{-(x|y)}\varphi _{{\bf s}+\varepsilon _j} (y)\Delta (y)^{-{N\over n}}m(dy) 
=\sum _{j=1}^n \alpha _j({\bf s})\Gamma _{\Omega } ({\bf s}+\varepsilon _j+\rho )
\varphi _{-{\bf s}-\varepsilon _j}(x),$$
one obtains 
\begin{eqnarray*}
{\rm tr}\, \bigl(\nabla \varphi _{-{\bf s}}(x)\bigr)
&=&-\sum _{j=1}^n \alpha _j({\bf s})
{\Gamma _{\Omega }({\bf s}+\varepsilon _j+\rho )\over \Gamma _{\Omega }({\bf s}+\rho )}
\varphi _{-{\bf s}-\varepsilon _j}(x) \\
&=&-\sum _{j=1}^n \alpha _j({\bf s})\Bigl(s_j-{d\over 4}(n-1)\Bigr)\varphi _{-{\bf s}-\varepsilon _j}(x),
\end{eqnarray*}
or
$${\rm tr}\, \bigl(\nabla \varphi _{\bf s}(x)\bigr)=\sum _{j=1}^n \alpha _j(-{\bf s})\Bigl(s_j+{d\over 4}(n-1)\Bigr)
\varphi _{{\bf s}-\varepsilon _j}(x).$$
In fact the explicit formula for $\Gamma _{\Omega}$,
$$\Gamma _{\Omega }({\bf s}+\rho )=(2\pi )^{N-n}\prod _{j=1}^n \Gamma \bigl(s_j-{d\over 4}(n-1)\bigr),$$
gives
$${\Gamma _{\Omega } ({\bf s}+\varepsilon _j+\rho )\over \Gamma _{\Omega }({\bf s}+\rho )}
={\Gamma \bigl(s_j+1-{d\over 4}(n-1)\bigr)\over \Gamma \bigl(s_j-{d\over 4}(n-1)\bigr)}
=s_j-{d\over 4}(n-1).$$

\bigskip

(ii) Recall that
$$D_{\nu }^{(2,-)}F(z)=-\langle z+e,\nabla F(z)\rangle -n\nu F(z).$$
From (i) we obtain
$$D_{\nu }^{(2,-)}\varphi _{\bf s}(z)
=\sum _{j=1}^n \bigl(s_j+{d\over 4}(n-1)\bigr)\alpha _j(-{\bf s})\varphi _{{\bf s}-\varepsilon _j}(z)
-\bigl(\sum _{j=1}^n s_j+n\nu \bigr)\varphi _{\bf s}(z).$$
By using $D_{\nu }^{(2,+)}=S_{\nu }^{(2)}D_{\nu }^{(2,-)}S_{\nu }^{(2)}$
and $S_{\nu }^{(2)}\varphi _{\bf s}(z)=\varphi _{-{\bf s}-\nu }(z),$
we get (ii).
\qed

We continue the proof of Theorem \ref{Pieri}. 
Let us write (ii) of Lemma \ref{TwoFormulas} with ${\bf s}={\bf k}-\rho $:
\begin{eqnarray*}
& &D_{\nu ,k}^{(2)}\Phi _{\bf k}(w)
=e^{i\theta }\Bigl(\sum _{j=1}^n (k_j+\nu -{d\over 2}(j-1)\bigr)\alpha _j({\bf k}-\rho )
\Phi _{{\bf k}+\varepsilon _j}(w)+|{\bf k}|\Phi _{\bf k}(w)\Bigr)\\
& &-e^{-i\theta }\Bigl(\sum _{j=1}^n \bigl(k_j+{d\over 2}(n-j)\bigr)\alpha _j(-{\bf k}+\rho )
\Phi _{{\bf k}-\varepsilon _j}(w)+(|{\bf k}|+n\nu )\Phi _{\bf k}(w)\Bigr).
\end{eqnarray*}
(Observe that $\sum _{j=1}^n\rho _j=0$.) Now, equaling the coefficient of $\Phi _{\bf k}(z)$ in both sides of (9.1), we obtain the formula of Theorem 8.1 
for all ${\bf s}={\bf m}+{\nu \over 2}-\rho $. Since both sides are polynomial functions in $\bf s$,
the equality holds for every~ $\bf s$. 


\bigskip

\noindent
{\it Acknowledgement.} The work of M.W. was partially supported by Grant-in-Aid for Challenging Exploratory Research No. 25610006 and by CREST, JST.



\bigskip

\begin{thebibliography}{99}

\bibitem{AAR1999}
Andrews, G. E., Askey, R. and Roy, R., 
\textit{Special functions}, 
Cambridge, 1999.

\bibitem{ADO2006}
Aristidou, M., Davidson, M. and \'Olafsson, G., 
\textit{Laguerre functions on symmetric cones and recursion relations in the real case},
J. Comput. Appl. Math. 199 (2006), 95--112.

\bibitem{BF1997}
Baker, T. H. and Forrester, P. J., 
\textit{The Calogero-Sutherland model and generalized classical polynomials},
Comm. Math. Phys. 188 (1997), 175--216.


\bibitem{BCKV2000}
Bump, D, Choi, K-K., Kurlberg, P. and Vaaler, J., 
\textit{A local Riemann hypothesis, I|}, 
Math. Z. 233 (2000), 1--19.

\bibitem{DO2003}
Davidson, M. and  \'Olafsson, G.,
\textit{Differential recursion relations for Laguerre functions
on Hermitian matrices|Integral transforms and special functions}, 
Integral Trans. Special Func., 14 (2003), 469--484.

\bibitem{DOZ2003}
Davidson, M., \'Olafsson, G. and Zhang, G., 
\textit{Laplace and Segal-Bargmann transforms on Hermitian symmetric spaces and orthogonal polynomials}, J. Func. Analysis 204 (2003), 157--195.

\bibitem{Di1990}
Dib, H., \textit{Fonctions de Bessel sur une alg\`ebre de Jordan}, 
J. Math. pures et appl. 69 (1990), 403--448. 

\bibitem{F2004}
Faraut, J., 
\textit{Analysis on the crown of a Riemannian symmetric space}, 
Amer. Math. Soc. Transl. 210 (2004), 99--110.

\bibitem{FK1994}
Faraut, J. and  Kor\'anyi, A., 
\textit{Analysis on symmetric cones}, Oxford 1994. 


\bibitem{H1963}
Hua, L.K., 
\textit{Harmonic analysis of functions of several variables in the classical 
domains}, Amer. Math. Soc. 1963.

\bibitem{L1998}
Lassalle, M., 
\textit{Coefficients binomiaux g\'en\'eralis\'es et polyn\^omes de Macdonald}, 
J. Funct. Analysis 158 (1998), 289--324. 

\bibitem{OZ1994}
\O rsted, B. and Zhang, G., \textit{Weyl quantization and tensor products of Fock and Bergman spaces}, 
Indiana Univ. Math. J. 43 (1994), 551--582.

\bibitem{OZ1995}
\O rsted, B. and  Zhang, G., 
\textit{Generalized principal series representations and tube domains}, 
Duke Math. J. 78 (1995), 335--357.

\bibitem{PZ1992}
Peetre, J. and Zhang, G., 
\textit{Aweighted Plancherel formula III. the case of a hyperbolic matrix domain}, 
Collect. Math. 43 (1992), 273--301.

\bibitem{SZ2007}
Sahi, S. and Zhang, G., \textit{Biorthogonal expansion of non-symmetric Jack functions}, 
SIGMA 3 (2007), 106, 9 pages.

\bibitem{S2000}
Schoutens, W., 
\textit{Stochastic Processes and Orthogonal Polynomials}, 
Lecture Notes in Statistics 146. Springer-Verlag, New York., 2000.

\bibitem{Z1995}
Zhang, G., \textit{Some recurrence formulas for spherical polynomials on tube domains}, 
Trans. Amer. Math. Soc. 347 (1995), 1725--1734.

\bibitem{Z2002}
Zhang, G., 
\textit{Invariant differential operators on symmetric cones and Hermitian symmetric spaces}, 
Acta Applicandae Mathematicae 73 (2002), 79--94.

\end{thebibliography}
\end{document}